\documentclass{birkjour}

\newtheorem{theorem}{Theorem}[section]
\newtheorem{proposition}[theorem]{Proposition}
\newtheorem{problem}[theorem]{Problem}

\newtheorem{remark}[theorem]{Remark}

\begin{document}

\title[Some notes on a method for proving inequalities by computer]
      {Some notes on a method for proving \\ inequalities by computer}

%----------Author 1
\author[Bojan Banjac]{Bojan Banjac}
\address{
         University of  Belgrade, Faculty of Electrical Engineering, \\
         Bulevar Kralja Aleksandra 73, 11000 Belgrade, Serbia;       \\[0.25 ex]
         University of Novi Sad, Faculty of Technical Sciences,      \\
         Trg Dositeja Obradovi\' ca 6, 21000 Novi Sad, Serbia
         }
\email{bojan.banjac@uns.ac.rs}
%----------Author 2
\author[Milica Makragi\' c]{Milica Makragi\' c}
\address{
         University of  Belgrade, Faculty of Electrical Engineering,\\
         Bulevar Kralja Aleksandra 73, 11000 Belgrade, Serbia
         }
\email{milica.makragic@etf.rs}
%----------Author 3
\author[Branko Male\v sevi\' c]{Branko Male\v sevi\' c}
\address{
         University of  Belgrade, Faculty of Electrical Engineering,\\
         Bulevar Kralja Aleksandra 73, 11000 Belgrade, Serbia
         }
\email{malesevic@etf.rs}
\thanks{Research is partially supported by the Ministry of Science and Education of the Republic of Serbia,
Grants ON 174032 and III 44006 }



\date{November 29, 2013} % \date{February 19, 2014}
%%% ----------------------------------------------------------------------

\begin{abstract}
In this article we consider mathematical fundamentals of one method for proving inequalities by computer,
based on the Remez algorithm. Using the well-known results of undecidability of the existence of zeros of
real elementary functions, we demonstrate that the considered method generally in practice becomes one heuristic for the verification of inequalities. We give some improvements of the inequalities considered in the theorems for which the existing proofs have been based on the numerical verifications of Remez algorithm.
\end{abstract}

%%% ----------------------------------------------------------------------
\maketitle
%%% ----------------------------------------------------------------------
%\tableofcontents

\medskip
$\!\!\!$
In this article we consider a method for proving inequalities by~computer, based on~the~Remez~algorithm,
which is presented in~the~article~\cite{Malesevic_07}.

\section{Mathematical fundamentals of Remezes based method \\ $\,\,\,\,\,\,\,\,$for proving inequalities}

\quad\enskip
Let $f\!:\![a,b] \longrightarrow R$ be a continuous function. The main topic of this article refers to inequalities
in the following form
\begin{equation}
\label{Ineq_Start}
f(x) \geq 0.
\end{equation}
In this section, we give mathematical fundamentals of a method for proving these inequalities \cite{Malesevic_07}.
Let us assume that there exist real numbers $n$ and $m$ such that there are finite and non-zero limits
\begin{equation}
\label{Limits}
\alpha = \lim\limits_{x \rightarrow a+}{\displaystyle\frac{f(x)}{(x-a)^n(b-x)^m}}
\quad \mbox{and} \quad
\beta  = \lim\limits_{x \rightarrow b-}{\displaystyle\frac{f(x)}{(x-a)^n(b-x)^m}}.
\end{equation}
The article \cite{Malesevic_07} considers the case when $n$ and $m$ are non-negative integers, where $n$ is the multiplicity of the root $x=a$ and $m$ is the multiplicity of the root $x=b$.
In the case wherein $x = a$ is not the root, then $n = 0$ and in the case wherein $x = b$ is not the root, then $m = 0$.
If for the function $f(x)$, at the point $x = a$, there is an approximation by Taylor polynomial of the $n$-th
degree and at the point $x = b$, there is an approximation
by Taylor polynomial of the $m$-th degree, then it follows that
\begin{equation}
\label{Limits_Taylor}
\alpha
=
\displaystyle\frac{f^{(n)}(a)}{n! \, (b-a)^{m}}
\quad \mbox{and} \quad
\beta
=
(-1)^m
\displaystyle\frac{f^{(m)}(b)}{m! \, (b-a)^{n}},
\end{equation}
which are determined in the article \cite{Malesevic_07}. Practically, for many examples of inequalities of the form
(\ref{Ineq_Start}), it is necessary to determine limits $\alpha$~and~$\beta$, for the function $f(x)$,
by formulas (\ref{Limits_Taylor}) for some non-negative integers $n$~and~$m$. Next, let us define
the following \mbox{$g$-function} for the real numbers $n, m$
and non-zero limits $\alpha, \beta$, as it was done in \cite{Malesevic_07}
\begin{equation}
\label{Function_g}
g(x)
=
g^{\,f}_{a,b}{(x)}
=
\left\{
\begin{array}{ccc}
\alpha                                   &\,\,:\,\,& x = a,                     \\[2.0 ex]
\displaystyle\frac{f(x)}{(x-a)^n(b-x)^m} &\,\,:\,\,& x \!\in\! (a,b),           \\[2.0 ex]
\beta                                    &\,\,:\,\,& x = b.
\end{array}
\right.
\end{equation}
The previous function is also continuous. Then, the following equivalence is true:
\begin{equation}
\label{Equivalence}
g(x) \geq 0 \;\Longleftrightarrow\; f(x) \geq 0.
\end{equation}
Thus, if $\alpha < 0$ or $\beta < 0$ the inequality (\ref{Ineq_Start}) is not true. Hence, we consider
only the cases $\alpha > 0$ and $\beta > 0$. Let us notice that if  the following implication is true:
\begin{equation}
\label{Implication}
g(x) > 0 \;\Longrightarrow\; f(x) \geq 0,
\end{equation}
then $f(x)$ can have roots at the end-points of the segment $[a,b]$. In next considerations we
use the function $g(x)$ instead of $f(x)$. One approach for proving the inequality $g(x) > 0$ is
based on the following statement.

\medskip

\begin{proposition}
\label{Prop_1}
Let $g:[a,b] \longrightarrow R$ be a continuous function. Then $g(x) > 0$ for $x \in [a,b]$ if and only if
there exist a polynomial $P(x)$ and a positive real number $\delta > 0$ such that$\,:$

\smallskip
\noindent
\begin{equation}
\label{Ineq_Appr}
| \, g(x) - P(x) \, | \leq \delta
\end{equation}
and

\medskip
\noindent
\begin{equation}
\label{Ineq_Stop}
P(x) - \delta > 0.
\end{equation}
\end{proposition}

\noindent
The necessary part of this statement is a simple consequence of the properties of the continuous functions
and the Weierstrass approximation theorem.

\medskip
\noindent
The polynomial $\hat{P}(x)$ of the $n$-th degree is {\it the minimax polynomial appro\-xi\-mation
of the $n$-th degree of the continuous function $g(x)$} over $[a,b]$, if the following
is true:
\begin{equation}
\max\limits_{x \in [a,b]}{| \, g(x) - \hat{P}(x) \, |}
\leq
\max\limits_{x \in [a,b]}{| \, g(x) - \tilde{P}(x) \, |},
\end{equation}
for arbitrary polynomial $\tilde{P}(x)$ of the $n$-th degree.
The minimax polynomial appro\-ximation $\hat{P}(x)$ of the $n$-th
degree is unique and the following well-known statement~is~true \mbox{\cite{J_F_Hart, Mayans_06}}.

\medskip

\begin{proposition}
$\![\,$The Chebyshev Equioscillation Theorem$\,]\!$
The polynomial $\hat{P}(x)$ is the minimax polynomial approximation of the $n$-th
degree of the continuous function $g(x)$ over~$[a,b]$ if and only if there exist~$n\!+\!2$ $t_{i}$-points
$a \!\leq\! t_{0} \!<\! \ldots \!<\! t_{n+1} \!\leq\! b$ such~that$\,:$
\begin{equation}
\label{de_la_Vallee-Poussin}
g(t_{i}) - \hat{P}(t_{i}) = (-1)^{i} \hat{\delta},
\quad
\hat{\delta} = \!\max\limits_{x \in [a,b]}{| \, g(x) - \hat{P}(x) \, |}.
\end{equation}
\end{proposition}

The polynomial Remez algorithm \cite{Remez_1934} is an iterative procedure. For the function
$g(x) : [a,b] \longrightarrow R$ and the selected accuracy $\varepsilon > 0$ the results are
the minimax polynomial appro\-ximation $\hat{P}(x)$ and the numerical estimate $\delta > 0$,
\cite{Fraser_Hart_62}, \cite{EW_Cheney_66}. All steps of Remez algorithm are described in details in
\cite{EW_Cheney_66, Robinson_05, Jean_Michel_Muller_06, Mastroianni_Milovanovic_08}.
Let us name {\em the step for extremes in the Remez algorithm} - searching for extremes
of the function $r(x) \!=\! P(x) \!-\! g(x)$.

\medskip
\noindent
The proof of effectiveness of the polynomial Remez algorithm is given in \cite{EP_Novodvorski_IS_Pinsker, EW_Cheney_66}.
In the case of twice continuously differentiable function $g(x)$, such that values $\pm\hat{\delta}$
from (\ref{de_la_Vallee-Poussin}) exist in the end-points $t_{0} \!=\! a$, $t_{n+1} \!=\! b$ and in $n$ interior
$t_{i}$-points also, the rate of convergence of the polynomial Remez algorithm is quadratic
\cite{Veidinger_60, J_F_Hart}. Based on the computer program {\sf Maple} and the {\sf numapprox} package,
the results of the polynomial Remez algorithm are, with the selected accuracy, the minimax polynomial approximation
$\hat{P}(x)$ and the numerical estimate of the absolute error $\hat{\delta}$ \cite{Geddes_93, Robinson_05}.
If it is not possible to determine the minimax polynomial approximation in the program, a message appears
that it is necessary to increase accuracy \cite{Geddes_93, B.Popov_O.Laushnyk, Robinson_05, Mayans_06}
and \cite{Malesevic_07}.

\smallskip
\noindent
Next, according to the Proposition \ref{Prop_1} for proving $g(x) \!>\! 0$ it is sufficient to use the polynomial
Remez algorithm of the minimax polynomial approximation of the function $g(x)$. On the other hand, it is possible
to use the varieties of rational Remez algorithms \cite{Cody_70, Litvinov_94}. In this case, we have
the well-known problems with convergence of these algorithms
\cite{PP_Petrushev_VI_Popov, Litvinov_94, Dunham_04}.

\smallskip
\noindent
Based on the previous considerations, we can determine a method from the article~\cite{Malesevic_07} more precisely by the following statement.

\medskip

\begin{proposition}
\label{Prop_2}
Let $f:[a,b] \longrightarrow R$ be the continuous~function~for~which~there
exist$\,:$

\smallskip
\noindent
{\boldmath $\,(i)\,\,$}
\begin{minipage}[t]{109 mm}
some real numbers $n$ and $m$, with limits $\alpha$ and $\beta$, determined~by~{\rm (\ref{Limits})},
as positive real numbers;
\end{minipage}

\smallskip
\noindent
{\boldmath $(ii)\,$}
\begin{minipage}[t]{109 mm}
the minimax polynomial approximation $P(x)$ and the numerical estimate of
the absolute~error~\mbox{$\delta \!>\! 0$} of the function $g(x) \!=\! g^{\,f}_{a,b}(x)$, determined by
{\rm (\ref{Function_g})}, such that {\rm (\ref{Ineq_Appr})} and {\rm (\ref{Ineq_Stop})} are true.
\end{minipage}

\smallskip
\noindent
Then, $f(x) \!\geq\! 0$ and $f(x)$ can have roots in the end-points of the~segment~$[a,b]$.
\end{proposition}
%
% \smallskip
%
% \begin{remark}
% Let us emphasize that the considered method, which is defined in the article {\rm \cite{Malesevic_07}}
% on the minimax polynomial approximations of function $g(x)$, can be modified for the other types of
% approximations of the function $g(x)$ in order to prove inequalities $g(x) \!>\! 0$.
% \end{remark}

\break

\section{An implementation of the method as a heuristic for computer verification of inequalities} % **********************************************

\medskip
\quad\enskip
In the article \cite{Malesevic_07}, the proofs of the considered inequalities in Theorems 1.2. and 1.10. are given by
the minimax linear approximations of the corresponding $g$-functions. The proofs of Theorems 1.2. and 1.10. from
\cite{Malesevic_07} are based on the assumptions that the appropriate minimax linear approximations and appropriate estimations
of the absolute errors, given by the program {\sf Maple}, are correctly given under the chosen accuracy.
Let us notice that, in practice, the correctness of the presented numerical proofs (by the presented method for the proving
inequalities \mbox{$f(x) \geq 0$} over $[a, b]$) is based on the correctness of all steps in realization of the
polynomial Remez algorithm, which is applied on the functions $g(x) \!=\! g^{\,f}_{a,b}{(x)}$ over~$[a,b]$.

\medskip
\noindent
Let us emphasize that in the article \cite{Malesevic_07} it has been remarked that the estimate of the absolute error
of Remez algorithm is of a numerical origin. This fact is also considered in the bibliography. P.$\,$L.~Richman
used the polynomial Remez algorithm in the article \cite{PL_Richman_69} and after deductions based on
the numerical estimate of the absolute error, he emphasized that for complicated functions, the numerical
estimate of the absolute error is not mathematically established bound (p.$\,$367.). The specified observation
of P.$\,$L.~Richman can be explained by the \mbox{following} statement.

\medskip

\begin{theorem}
The step for extremes in the Remez algorithm is reduced to an undecidable problem, if in that step,
searching for extremes of the function $r(x) \!=\! P(x) \!-\! g(x)$ is determined by zeros of
the first derivative.
\end{theorem}

\smallskip
\noindent
{\bf Proof.}
It follows from the well-known result of P.$\,$S. Wang of the undecidability of the existence of zeros of
real elementary functions \cite{Wang_74} {\big (}\hspace*{0.25 mm}see Section$\;$9 in~\cite{Poonen_12}
and~\cite{Encyclopedia_of_Mathematics}\hspace*{0.25 mm}{\big )}. $\Box$

\medskip

\begin{remark}
Let us to be emphasize, that in various numerical algorithms generally there appears undecidability of those steps in which
we seek zeros of complicated functions.
\end{remark}

\medskip

\noindent
In general, let us notice that the class of continuous functions of one variable, over some segment, for
all root-finding algorithms, is incomplete.
% Then for this class functions we considered root-finding algorithms as heuristics
% in sense rules that may find a root for functions over this class but are not guaranteed to.
Hence, in practice, for complicated functions the presented method for proving inequalities
becomes an inequalities verification heuristic. Let us emphasize that
in the articles \cite{Malesevic_07}, \cite{Chen_Cheung_Wang_2011}, \cite{Chen_Cheung_2011} proofs are provided, based on calculations done by the program {\sf Maple} without proofs
of correctness of the numerical estimates of appropriate absolute errors. Based on the previous
facts on implementation of the considered method in~this~section, we give new proofs of
Theorem~1.2 from~\cite{Malesevic_07}, Theorem~3.2 from~\cite{Chen_Cheung_Wang_2011} and Theorem~2
from~\cite{Chen_Cheung_2011} instead of the existing numerical verifications. In this article
Theorems \ref{TH1} and \ref{TH3} provide some new inequalities.

\smallskip

\begin{theorem}\label{TH1}
Let $K(x) \!=\! \!\!\mbox{\small $\displaystyle\int\limits_{0}^{\infty}{\!\!e^{-t}
\displaystyle\frac{t^{x}\!-\!1}{t\!-\!1} \: dt}$} \; (x \!\geq\! 0)$
be the Kurepa's \mbox{function}~\mbox{\rm \cite{Kurepa_73}},~{\rm \cite{Ivic_Mijajlovic_95, Malesevic_10}};
then, for values $x \!\in\! [0,1]$, the following inequality is true$:$
\begin{equation}
\label{Ineq_1}
K(x) \leq \displaystyle K'\!(0) \, x  \, + {\big (}1 - K'(0){\big )}x^2,
\end{equation}
where $K'(0) \!=\! 1.432 \, 205 \, ... \,$.
\end{theorem}

\smallskip

\noindent
{\bf Proof.} Let $H_{2}(x) = \alpha + \beta x + \gamma x^2$ be the Hermite polynomial at the
nodes $x=0$ (with multiplicity two) and $x=1$, of the function $K(x)$ with
\begin{equation}
\alpha = K(0) = 0,
\quad
\beta = K'(0),
\quad
\alpha + \beta + \gamma = K(1) = 1.
\end{equation}
Then, for $x \!\in\! [0,1]$ we have
\begin{equation}
K(x) - H_{2}(x)
=
\displaystyle\frac{K'''(\xi)}{3!}\,x^2(x-1),
\quad
0 < \xi < 1,
\end{equation}
where $K'''(x) \!=\! \mbox{\small $\displaystyle\int\limits_{0}^{\infty}{\!\!e^{-t}
\displaystyle\frac{t^{x} \log^3t}{t\!-\!1} \: dt}$} > 0$. Therefore, we obtained
inequality~(\ref{Ineq_1}).$\,\Box$

\begin{remark}
Let us emphasize that the inequality {\rm (\ref{Ineq_1})} is an improvement of the inequality
$(1.8)$ from Theorem $1.2.$ of the article {\rm \cite{Malesevic_07}} $($see also \cite{Malesevic_04}$)$.
\end{remark}

\smallskip
\noindent
Let us remark that Theorem 1.10. from the article \cite{Malesevic_07} is proved in the article
\cite{MalesevicB_07} and is considered in articles
\cite{Zhu_07a, Zhu_07, Zhu_08, Zhu_09, Qi_Guo_12, Qi_Luo_Guo_12, Guo_Luo_Qi_2013}.
Inequalities of similar type, such as an inequality which is considered in Theorem 1.10~\cite{MalesevicB_07},
for various trigonometric functions, have been considered in the articles \cite{MitrinovicVasic_70, Qi_Niu_Guo_2009, Qi_Zhang_Guo_2009, Zhang_Guo_2009, Guo_Qi_2010, Sun_Zhu_11, Chen_Cheung_2012, Zhao_Wei_Guo_Qi_2012, Chen_Sandor_13, Deng_Chen_2014, Debnath_Zhu_Mortici_2014}. An interesting application of this type inequalities can be seen in \cite{Malesevic_97},
and is considered~in~the~article~\cite{G.T.F.de_Abreu_09}.

\begin{theorem}\label{TH2}
\hspace*{-1.5 mm} {\rm \cite{Chen_Cheung_Wang_2011}}$\,$ For $0 \leq x \leq 1$,
\begin{equation}
\label{CCW_inequality}
\displaystyle\frac{(\pi/2)(1-x)^{(\pi+2)/\pi^2}}{(1+x)^{(\pi-2)/\pi^2}} \leq \mbox{\rm arccos} \; x \, .
\end{equation}
\end{theorem}

\noindent
{\bf Proof.}
Let us notice that the proof of the previous inequality in \cite{Chen_Cheung_Wang_2011} is based
on the following inequality
\begin{equation}
\label{CCW_g}
g(t)=4 t \cos (2\,t) + 2 \pi \, t - \frac{\pi^2}{2} \sin (2\,t) \geq 0, \qquad 0 \leq t \leq \frac{\pi}{4}.
\end{equation}
We provide a new proof of the previous inequality. Let us denote $u = 2\,t$, then the inequality becomes

\vspace*{-2.5 mm}

\begin{equation}
\varphi(u) = g(2\,t) = 2 u \cos (u) + \pi \, u - \frac{\pi^2}{2} \sin (u) \geq 0, \qquad 0 \leq u \leq \frac{\pi}{2}.
\end{equation}
Therefore,
$
\varphi(u) \!\geq\! 0
\Longleftrightarrow
2 u \cos u + \pi \, u
\!\geq\!
\mbox{\small $\displaystyle\frac{\pi^2}{2}$} \sin u
\Longleftrightarrow
u \geq \displaystyle\frac{\mbox{\small $\displaystyle\frac{\pi}{2}$} \sin u}{
1 + \mbox{\small $\displaystyle\frac{2}{\pi}$} \cos u},
$
for $0 \!\leq\! u \!\leq\! \mbox{\small $\displaystyle\frac{\pi}{2}$}$, which is true
according~to \cite{Zhu_07} (Theorem 7.), see also \cite{Qi_Luo_Guo_12}.$\;\;\Box$

\begin{remark}
Let us emphasize that the previous proof is an improvement of the proof of the inequality $(3.6)$
from Theorem $3.2.$ of the article {\rm \cite{Chen_Cheung_Wang_2011}}.
\end{remark}

\break

\medskip

\begin{theorem}\label{TH3}
For $0 < x < \pi/2$, the following inequalities are true:
\begin{equation}
\label{Ineq_21}
\left(\frac{\pi^2}{\pi^2-4x^2}\right)^{\!\alpha}
\!<
\displaystyle
\frac{
\left(
\mbox{\small $\displaystyle\frac{{\pi}^{2}}{18}$}
-
\mbox{\small $\displaystyle\frac{2}{3}$}
\right){x}^{4}
+
\left(
\mbox{\small $\displaystyle\frac{{\pi}^{2}}{3}$}
-
4
\right){x}^{2}
+
{\pi}^{2}}{
{\pi }^{2}-4\,{x}^{2}}
\,<\,
\displaystyle
\frac{\tan x}{x}
\end{equation}

\noindent
and

\vspace*{-4.5 mm}

\begin{equation}
\label{Ineq_22}
\frac{\tan x}{x}
<
\displaystyle
\frac{
-\mbox{\small $\displaystyle\frac{1}{{\pi}^{2}}$}\,x^{6}
+\mbox{\small $\displaystyle\frac{1}{2}$}\,x^{4}
-\mbox{\small $\displaystyle\frac{{\pi}^{2}}{16}$}\,x^{2}
+{\pi}^{2}}{
{\pi }^{2}-4\,{x}^{2}}
<
\left(\frac{\pi^2}{\pi^2-4x^2}\right)^{\!\beta}
\end{equation}
with the constants $\alpha = \mbox{\small $\displaystyle\frac{\pi^2}{12}$} = 0.822 \, 467 \, \ldots$ and $\beta=1$ .
\end{theorem}

\noindent
{\bf Proof.} Let us denote:

\vspace*{-5.0 mm}

\begin{equation}
\begin{array}{rcl}
f_1(x)
\!&\!\!=\!\!&\!
\left(\mbox{\small $\displaystyle\frac{\pi^2}{\pi^2-4x^2}$}\right)^{\!\!\frac{\pi^2}{12}}       \\[2.5 ex]
f_2(x)
\!&\!\!=\!\!&\!
\displaystyle
\frac{
\left(
\mbox{\scriptsize $\displaystyle\frac{{\pi}^{2}}{18}$}
-
\mbox{\scriptsize $\displaystyle\frac{2}{3}$}
\right)\mbox{\small ${x}^{4}$}
+
\left(
\mbox{\scriptsize $\displaystyle\frac{{\pi}^{2}}{3}$}
-
\mbox{\small $4$}
\right)\mbox{\small ${x}^{2}$}
+
\mbox{\small ${\pi}^{2}$}}{
\mbox{\small ${\pi }^{2}$}-\mbox{\small $4\,{x}^{2}$}}                                            \\[3.0 ex]
f_{3}(x)
\displaystyle
\!&\!\!=\!\!&\!
\mbox{\small $\displaystyle\frac{\tan x}{x}$}                                                     \\[1.0 ex]
f_{4}(x)
\displaystyle
\!&\!\!=\!\!&\!
\displaystyle
\frac{
-\mbox{\scriptsize $\displaystyle\frac{1}{{\pi}^{2}}$}\,\mbox{\small $x^6$}
+\mbox{\scriptsize $\displaystyle\frac{1}{2}$}\,\mbox{\small $x^4$}
-\mbox{\scriptsize $\displaystyle\frac{{\pi}^{2}}{16}$}\,\mbox{\small $x^2$}
+\mbox{\small ${\pi}^{2}$}}{
\mbox{\small ${\pi }^{2}$}-\mbox{\small $4\,{x}^{2}$}}                                            \\[1.0 ex]
f_{5}(x)
\!&\!\!=\!\!&\!
\displaystyle
\mbox{\small $\displaystyle\frac{\pi^2}{\pi^2-4x^2}$}
\end{array}
\end{equation}
All details of proofs of the next four inequalities $f_{1}(x) \!<\! f_{2}(x) \!<\! f_{3}(x) \!<\! f_{4}(x) \!<\! f_{5}(x)$
for \mbox{$0 \!< \!x \!<\! \mbox{\small $\displaystyle\frac{\pi}{2}$}$ }, are
presented in the Appendix of this article. The following conclusions are true:

\smallskip
\noindent
{\boldmath $1^{0}$}
Inequality $f_{2}(x) \!>\! f_{1}(x)$ is true on the following facts:

\smallskip
\noindent
$
\;\;\;
\begin{array}{ll}
\mbox{\bf \boldmath $a$)} \!\!&\!\! \mbox{$f_{1}(x)$ and $f_{2}(x)$ are convex functions;}               \\[1.5 ex]
\mbox{\bf \boldmath $b$)} \!\!&\!\! \lim\limits_{x \rightarrow 0+}{\!{\big (}f_2(x)\!-\!f_1(x){\big )}} \!=\! 0;   \\[1.5 ex]
\mbox{\bf \boldmath $c$)} \!\!&\!\! \lim\limits_{x \rightarrow 0+}{\!{\big (}f_2^{(k)}(x)\!-\!f_1^{(k)}(x){\big )}}
\!=\! 0,  \, \mbox{for} \; k=1..5;                                                                                 \\[1.0 ex]
\mbox{\bf \boldmath $d$)} \!\!&\!\! \lim\limits_{x \rightarrow 0+}{\!{\big (}f_2^{(6)}(x)\!-\!f_1^{(6)}(x){\big )}}
\!=\! \mbox{\small $\displaystyle\frac{40(144-\pi^4)}{9 \pi^4}$} \!=\!  2.125 ... > 0;                             \\[1.5 ex]
\mbox{\bf \boldmath $e$)} \!\!&\!\! \lim\limits_{x \rightarrow \frac{\pi}{2}-0}{\!f_1(x)} \!=\! +\infty
\;\;\;\mbox{and}\;
\lim\limits_{x \rightarrow \frac{\pi}{2}-0}{\!{\big (}f_2(x)\!-\!f_1(x){\big )}} \!=\! +\infty.
\end{array}
$

\smallskip
\noindent
{\boldmath $2^{0}$}
Inequality $f_{3}(x) \!>\! f_{2}(x)$ is true, according to a method of proving a class of trigonometric inequalities
based on approximations of the sine and cosine functions by Maclaurin polynomials \cite{Mortici_2011} and
\cite{Malesevic_Makragic_2015}, see the Appendix.

\smallskip
\noindent
{\boldmath $3^{0}$}
Inequality $f_{4}(x) \!>\! f_{3}(x)$ is true, according to a method of proving a class of trigonometric inequalities
based on approximations of the sine and cosine functions by Maclaurin polynomials \cite{Mortici_2011} aand
\cite{Malesevic_Makragic_2015}, see the Appendix.

\medskip
\noindent
{\boldmath $4^{0}$}
Inequality $f_{5}(x) \!>\! f_{4}(x)$ is trivial, see the Appendix. $\Box$

\begin{remark}
Let us emphasize that the inequalities {\rm (\ref{Ineq_21})} and {\rm (\ref{Ineq_22})} are an improvement
of the inequality $(4)$ from Theorem $2.$ of the article {\rm \cite{Chen_Cheung_2011}}.
\end{remark}

\break

\medskip
At the end of this article, let us emphasize that the methods of the high-accuracy computations
of the minimax polynomial approximations and the appropriate estimations of the absolute errors
are considered in Chapter~$5$ of the Book \cite{F_Bornemann_D_Laurie_S_Wagon_J_Waldvogel}.
Now, with reference to the Problem 5 of the SIAM hundred$\,\!$--$\!\,$digit challenge
\cite{F_Bornemann_D_Laurie_S_Wagon_J_Waldvogel}, we state this as an open problem.

\begin{problem}
Let $f(z) = 1/\Gamma(z)$ where $\Gamma(z)$ is the gamma function, and let $\hat{p}(z)$ be the
cubic polynomial that best approximates $f(z)$ on the unit disk in the supremum norm
$\|\,.\,\|_{\infty}.\,$Prove correctness of all steps in computing
$
\|f(z) \!\,-\,\! \hat{p}(z)\|_{\infty} = 0.214 \, 335 \!\; \ldots
$
from~{\rm \cite{F_Bornemann_D_Laurie_S_Wagon_J_Waldvogel}}.
\end{problem}

\smallskip
{\sc Acknowledgement.} The authors would like to thank anonymous reviewer for his/her valuable comments
and suggestions, which were helpful in improving the article.

\break

$\,$

\bigskip
\noindent
\hspace*{4.5 mm}
{\large \bf The Appendix}

\bigskip

\bigskip
\noindent
{\boldmath $1^{0}$}
The conclusions {\bf \boldmath $a$)} -- {\bf \boldmath $e$)} follow from the facts:

\bigskip
\noindent
{\boldmath $1^{0}/\mbox{\small \bf \boldmath $a$)}$}
For $0 < x <\mbox{\small $\displaystyle\frac{\pi}{2}$}$ we have
\begin{equation}
f_{2}''(x)
=
\displaystyle\frac{
2{\Big (}
{\big (}8\pi^2\!-\!96{\big )}x^6 +
{\big (}-6\pi^4\!+\!72\pi^2{\big )}x^4+
3\pi^6x^2+
3\pi^6
{\Big )}}{
{\big (}\pi^2\!-\!4x^2{\big )}^3}
>
0
\end{equation}
and
\begin{equation}
f_{1}''(x)
=
\displaystyle\frac{
2\left(
{\Big (}
2\pi^{\left(\frac{\pi^2}{6}+4\right)}+12\pi^{\left(\frac{\pi^2}{6}+2\right)}
{\Big )} x^2
+
3\pi^{\left(\frac{\pi^2}{6}+4\right)}
\right)}{
9{\big (}\pi^2-4x^2{\big )}^{2-\frac{\pi^2}{12}}}
>
0.
\end{equation}
\bigskip
\noindent
{\boldmath $1^{0}/\mbox{\small \bf \boldmath $b$)}\mbox{\rm ,}\,\mbox{\small \bf \boldmath $c$)}\mbox{\rm ,}\,\mbox{\small \bf \boldmath $d$)}$}
For $x \rightarrow 0_{+}$ we have
\begin{equation}
f_2(x)
=
1+\frac{1}{3}x^2+\frac{\pi^2\!+\!12}{18\pi^2}x^4+\frac{2(\pi^2\!+\!12)}{9\pi^4}x^6+o(x^6)
\end{equation}
and
\begin{equation}
f_1(x)
=
1+\frac{1}{3}x^2+\frac{\pi^2\!+\!12}{18\pi^2}x^4+\frac{\pi^4\!+\!36\pi^2\!+\!288}{162\pi^4}x^6+o(x^6),
\end{equation}
with non-zero values of the appropriate derivatives
$f_{2}^{(2)}(0)\!=\!f_{1}^{(2)}(0)\!=\!\mbox{\footnotesize $\displaystyle\frac{2}{3}$}$,
$f_{2}^{(4)}(0)\!=\!f_{1}^{(4)}(0)\!=\!\mbox{\footnotesize $\displaystyle\frac{4(\pi^2\!+\!12)}{3\pi^2}$}$,
$
f_{2}^{(6)}(0)
\!=\!
\mbox{\footnotesize $\displaystyle\frac{160(\pi^2\!+\!12)}{\pi^4}$}
$,
$
f_{1}^{(6)}(0)
\!=\!
\mbox{\footnotesize $\displaystyle\frac{40(\pi^4\!+\!36\pi^2\!+\!288)}{9\pi^4}\,$}.
$

\bigskip
\noindent
{\boldmath $1^{0}/\mbox{\small \bf \boldmath $e$)}$} For $x \rightarrow \mbox{\small $\displaystyle\frac{\pi}{2}$}_{-}$,
it is true
\begin{equation}
f_{2}(x)
\sim
\displaystyle\frac{\left(\pi/2\right)^{\!\frac{\;\pi^2}{12}}}{\left(\pi-2x\right)^{\!\frac{\;\pi^2}{12}}}
\quad\mbox{and}\quad
f_{1}(x)
\sim
\displaystyle\frac{\mbox{\scriptsize $\frac{\pi^3(\pi^2\!+12)}{576}$}}{\pi-2x}\,.
\end{equation}

\bigskip
\noindent
{\boldmath $2^{0}$}
In this section of this appendix we give a proof based on the method from \cite{Malesevic_Makragic_2015} for the following equivalent inequality:
\begin{equation}
\label{Ineq_23}
\varphi(x)\!=\!\big{(}\!\underbrace{\pi^{2}-4x^{2}}_{>0}\!\big{)}\sin \!x+
\Bigg{(}\!\!\!\underbrace{\Big{(}\!\mbox{\scriptsize $\displaystyle\frac{\pi^{2}}{18}-\frac{2}{3}$}
\!\Big{)}x^{4}+\Big{(}\!\mbox{\scriptsize $\displaystyle\frac{\pi^{2}}{3}-4$}\!\Big{)}x^{2}
+\pi^{2}\!\!\Bigg{)}x}_{>0}
(\!-\!\cos \!x\!)>0
\end{equation}
for $x \in (0,\frac{\pi}{2})$.

\smallskip
\noindent
Let us consider two cases:

\medskip
\noindent
{\bf \boldmath $a$)}
If $x \!\in \!(0,\frac{\pi}{4}]$ the following holds: $\sin x \geq \underline{T}_{\,7}^{\sin,0}(x)$ and
$-\cos x \geq -\overline{T}_{4}^{\,\cos,0}(x)$, \cite{Malesevic_Makragic_2015}.
Then for $x \in (0,\frac{\pi}{4}]$ it is valid:
\begin{equation}
\label{Ineq_24}
\begin{array}{rcl}
\varphi(x)
\!&\!\!\!>\!\!\!&\!
\!\Big{(}\!\underbrace{\pi^{2}\!-\!4x^{2}}_{>0}\!
\Big{)}\underline{T}_{\,7}^{\sin,0}(x)
                                                                                             \\[1.5ex]
\!&\!\!\!+\!\!\!&\!
\Bigg{(}\!\!\!\underbrace{\Big{(}\!\mbox{\scriptsize $\displaystyle\frac{\pi^{2}}{18}\!-\!\frac{2}{3}$}
\!\Big{)}x^{4}\!+\!\Big{(}\!\mbox{\scriptsize $\displaystyle\frac{\pi^{2}}{3}\!-\!4$}\!\Big{)}x^{2}
\!+\!\pi^{2}\!\!\Bigg{)}x}_{>0}
\big{(}\!\!-\!\overline{T}_{4}^{\,\cos,0}(x)\!\big{)}\!=\!P_{9}(x),
\end{array}
\end{equation}
where $P_{9}(x)$ is the polynomial
\begin{equation}
\label{Ineq_25}
\begin{array}{rcl}
P_{9}(x)
\!&\!\!=\!\!&\!
\mbox{\footnotesize $-\displaystyle\frac{x^{5}}{15120}$}
{\big (}
\mbox{\small $(35\pi^{2}\!-\!432)x^{4}\!+\!(3024\!-\!207\pi^{2})x^{2}\!+\!10080\!-\!1176\pi^{2}$}
{\big )}
                                                                                  \\[2.0 ex]
\!&\!\!=\!\!&\!
\mbox{\footnotesize $-\displaystyle\frac{x^{5}}{15120}$}
P_{4}(x).
\end{array}
\end{equation}
Then we determine the sign of the polynomial $P_{4}(x)$ for $x \in (0,\frac{\pi}{4}]$.
A real numerical factorization of the polynomial $P_{4}(x)$, has been determined via {\sf Matlab}
software, and given with
$P_{4}(x)\!=\!\alpha(x\!-\!x_{1})(x\!-\!x_{2})(x\!-\!x_{3})(x\!-\!x_{4})$,
where
$
\alpha\!\!=\!\!-86.563\ldots,
x_{1}\!\!=\!\!1.364\ldots,\!
$
$
x_{2}\!\!=\!\!3.077\ldots,\!
$
$
x_{3}\!\!=\!\!-x_{1}\!
$
and
$
x_{4}\!\!=\!\!-x_{2}.
$
The polynomial $P_{4}(x)$ has exactly four simple real roots with a symbolic radical representation and
the corresponding numerical values $x_{1}$, $x_{2}$, $x_{3}$ and $x_{4}$.

\noindent
Since $P_{4}(0)<0$ it follows that $P_{4}(x)<0$ for $x \in (0,\frac{\pi}{4}]\subset(x_{3},x_{1})$.
Finally, we conclude that
$P_{9}(x) \!>\! 0 \,\,\mbox{\rm for}\,\, x \!\in\! (0,\!\frac{\pi}{4}]
\Longrightarrow \varphi(x) \!>\! 0\,\,\mbox{\rm for}\,\, x \!\in\! (0,\frac{\pi}{4}]$.

\smallskip
\noindent
{\bf \boldmath $b$)}
If $x \!\in \!(\frac{\pi}{4},\frac{\pi}{2})$, let us define the function
\begin{equation}
\label{Ineq_26}
\begin{array}{rcl}
\phi(x)
          \!&\!\!\!\!=\!\!\!\!&\!
\varphi(\frac{\pi}{2}-x)=
(4\pi x-4x^{2})\cos x
                                                                                          \\[1.5 ex]
           \!&\!\!\!\!+\!\!\!\!&\!
(\frac{\pi}{2}-x)\Big{(}\!\big{(}\frac{\pi^{2}}{18}\!-\!\frac{2}{3}\!\big{)}x^{4}\!+\!
\big{(}\!\frac{4\pi}{3}\!-\!\frac{\pi^{3}}{9}\!\big{)}x^{3}\!+\!
\big{(}\!\frac{\pi^{4}}{12}\!-\!\frac{2\pi^{2}}{3}\!-\!4\!\big{)}x^{2}\!+\!
\big{(}\!4\pi\!-\!\frac{\pi^{5}}{36}\!\big{)}x
                                                                                            \\[1.5 ex]
           \!&\!\!\!\!+\!\!\!\!&\!
\frac{\pi^{6}}{288}\!+\!\frac{\pi^{4}}{24}\Big{)}(\!-\!\sin x).
\end{array}
\end{equation}
Now we prove that $\phi(x)\!>\!0$ for $x \!\in \!(0,\frac{\pi}{4})$. The following holds:
$\cos x\!\geq \!\underline{T}_{\,6}^{\cos,0}(x)$ and $-\sin x\geq-\overline{T}_{5}^{\,\sin,0}(x)$, \cite{Malesevic_Makragic_2015}.
Then for $x\in \left(0,\frac{\pi}{4}\right)$ it holds:
\begin{equation}
\label{Ineq_27}
\begin{array}{rcl}
\phi(x)
       \!&\!\!\!\!>\!\!\!\!&\!
(4\pi x-4x^{2})\underline{T}_{\,6}^{\cos,0}(x)
                                                                                          \\[1.5 ex]
           \!&\!\!\!\!+\!\!\!\!&\!
(\frac{\pi}{2}-x)\Big{(}\!\big{(}\frac{\pi^{2}}{18}\!-\!\frac{2}{3}\!\big{)}x^{4}\!+\!
\big{(}\!\frac{4\pi}{3}\!-\!\frac{\pi^{3}}{9}\!\big{)}x^{3}\!+\!
\big{(}\!\frac{\pi^{4}}{12}\!-\!\frac{2\pi^{2}}{3}\!-\!4\!\big{)}x^{2}\!+\!
\big{(}\!4\pi\!-\!\frac{\pi^{5}}{36}\!\big{)}x
                                                                                            \\[1.5 ex]
           \!&\!\!\!\!+\!\!\!\!&\!
\frac{\pi^{6}}{288}\!+\!\frac{\pi^{4}}{24}\Big{)}(-\overline{T}_{5}^{\,\sin,0}(x))\!=\!Q_{10}(x),
\end{array}
\end{equation}
where $Q_{10}(x)$ is the polynomial
$$
\begin{array}{rcl}
Q_{10}(x)
           \!&\!\!\!\!=\!\!\!\!&\!
\mbox{\footnotesize $\displaystyle-\frac{x}{69120}$}
{\Big (}\!\!
\left( 384-32\,{\pi }^{2} \right) {x}^{9}+ \left( 80\,{\pi }^{3}-960\,\pi  \right) {x}^{8}
                                                                                               \\[2.5 ex]
            \!&\!\!\!\!+\!\!\!\!&\!
 \left( \!1408\,{\pi }^{2}\!-\!5760-80\,{\pi }^{4} \right) {x}^{7}\!+\!
 \left( \!40\,{\pi }^{5}\!+\!16128\,\pi -1792\,{\pi }^{3}\! \right) {x}^{6}
                                                                                               \\[2.5 ex]
             \!&\!\!\!\!+\!\!\!\!&\!
\left( -10\,{\pi }^{6}+11520-18048\,{\pi }^{2}+1576
\,{\pi }^{4} \right) {x}^{5}
\end{array}
$$

\break

\begin{equation}
\label{Ineq_28}
\begin{array}{rcl}
             \!&\!\!\!\!+\!\!\!\!&\!
\left( 13440\,{\pi }^{3}-788\,{\pi }^{5}
-57600\,\pi +{\pi }^{7} \right) {x}^{4}
                                                                                                \\[1.5 ex]
              \!&\!\!\!\!+\!\!\!\!&\!
\left( 69120\,{\pi }^{2}+
138240-9120\,{\pi }^{4}+200\,{\pi }^{6} \right) {x}^{3}
                                                                                                \\[1.5 ex]
              \!&\!\!\!\!+\!\!\!\!&\!
\left( -23040
\,{\pi }^{3}-276480\,\pi -20\,{\pi }^{7}+4560\,{\pi }^{5} \right) {x}^
{2}+
                                                                                                \\[1.5 ex]
              \!&\!\!\!\!+\!\!\!\!&\!
\left( 276480+138240\,{\pi }^{2}-2880\,{\pi }^{4}-1200\,{\pi }^{6
} \right) x+1440\,{\pi }^{5}
                                                                                                \\[1.5 ex]
              \!&\!\!\!\!+\!\!\!\!&\!
120\,{\pi }^{7}-276480\,\pi
{\Big )}
=\mbox{\footnotesize $\displaystyle-\frac{x}{69120}$}Q_{9}(x).
\end{array}
\end{equation}
Then we determine the sign of the polynomial
$Q_{9}(x)$ for $x\in \left(0,\frac{\pi}{4}\right)$.
Let us look at the fifth derivative of the polynomial $Q_{9}(x)$, as the fourth degree polynomial, in the following form:
\begin{equation}
\begin{array}{rcl}
\label{Ineq_29}
Q^{(5)}_{9}(x)
       \!&\!\!\!\!=\!\!\!\!&\!
\left( -483840\,{\pi }^{2}\!+\!5806080 \right) {x}^{4}\!+\!
\left( -6451200\,\pi \!+\!537600\,{\pi }^{3} \right) {x}^{3}
                                                                                                \\[1.5 ex]
       \!&\!\!\!\!+\!\!\!\!&\!
\left( -201600\,{\pi }^{4}\!+\!3548160\,{\pi }^{2}\!-\!14515200 \right) {x}^{2}
                                                                                                \\[1.5 ex]
        \!&\!\!\!\!+\!\!\!\!&\!
\left( 11612160\,\pi \!+\!28800\,{\pi }^{5}\!-\!1290240\,{\pi }^{3} \right) x
-2165760\,{\pi }^{2}
                                                                                                \\[1.5 ex]
         \!&\!\!\!\!-\!\!\!\!&\!
1200\,{\pi }^{6}+1382400+189120\,{\pi }^{4}
\,.
\end{array}
\end{equation}
A real numerical factorization of the polynomial $Q^{(5)}_{9}(x)$, has been determined via {\sf Matlab} software, and given with
$Q^{(5)}_{9}(x)
\!=\!
\beta(x\!-\!x_{1})(x\!-\!x_{2})(x\!-\!x_{3})(x\!-\!x_{4})$,\\
where
$
\beta\!=\!1030770.606\ldots,
x_{1}\!=\!0.566\ldots,\,\,
x_{2}\!=\!1.589\ldots,\,\,
$
$
x_{3}\!=\!2.506\ldots,
x_{4}=\!-1.171\ldots.
$
The polynomial equation $Q^{(5)}_{9}(x)=0$ has got exactly four simple real roots with a symbolic radical representation and the corresponding numerical values $x_{1}$, $x_{2}$, $x_{3}$ and $x_{4}$.
Let us notice that $x_{1}=0.566\ldots \in (0,\frac{\pi}{4})$.
Since
$Q^{(5)}_{9}(0)=-2165760\,{\pi }^{2}-1200\,{\pi }^{6}+1382400+189120\,{\pi }^{4}<0$ and $Q^{(5)}_{9}(1)=-7326720-752640\pi^{3}+898560\pi^{2}-1200\pi^{6}+28800\pi^{5}-12480\pi^{4}+5160960\pi>0$
it follows that $Q^{(5)}_{9}(x)<0$ for $x<x_{1}$ and $Q^{(5)}_{9}(x)>0$ for $x \in (x_{1},x_{2})$.
Therefore, $Q^{(4)}_{9}(x)$ is a monotonically decreasing function for  $x \in (0,x_{1})$
and a monotonically increasing function for $x \in (x_{1},\frac{\pi}{4}) \subset (x_{1},x_{2})$.
Hence, $Q^{(4)}_{9}(x)$ reaches the minimum at the point $x_{1}\!=\!0.566\ldots$.
Then, since $Q^{(4)}_{9}(0)=24\pi^{7}-1382400\pi-18912\pi^{5}+322560\pi^{3}<0$ and
$Q^{(4)}_{9}(\frac{\pi}{4})=1362\pi^{5}+68400\pi^{3}-1036800\pi+\frac{9}{2}\pi^{7}<0$,
it follows that $Q^{(4)}_{9}(x)$ is a negative function for $x \in (0,\frac{\pi}{4})$. Therefore, $Q^{(3)}_{9}(x)$ is a monotonically decreasing function for $x \in (0,\frac{\pi}{4})$.
Then, since $Q^{(3)}_{9}(\frac{\pi}{4})>0$, it follows that $Q^{(3)}_{9}(x)$ is a positive function for
$x \in (0,\frac{\pi}{4})$, so it follows that $Q^{(2)}_{9}(x)$ is a monotonically increasing function for
$x \in (0,\frac{\pi}{4})$. Since $Q^{(2)}_{9}(\frac{\pi}{4})<0$, then it holds that $Q^{(2)}_{9}(x)$ is a negative function for $x \in (0,\frac{\pi}{4})$. Hence, it follows that $Q'_{9}(x)$ is a monotonically decreasing function for $x \in (0,\frac{\pi}{4})$.
Since $Q'_{9}(\frac{\pi}{4})>0$, then it is valid that $Q'_{9}(x)$ is a positive function for
$x \in (0,\frac{\pi}{4})$. Hence, it holds that $Q_{9}(x)$ is a monotonically increasing function
for $x \in (0,\frac{\pi}{4})$. Finally, since
$Q_{9}(\frac{\pi}{4})<0$ we conclude that
$Q_{9}(x) \!<\! 0  \,\,\mbox{\rm for}\,\, x \!\in\! (0,\frac{\pi}{4})
\Longrightarrow
Q_{10}(x) \!>\! 0 \,\,\mbox{\rm for}\,\, x \!\in\! (0,\frac{\pi}{4})
\Longrightarrow
\phi(x)\!>\! 0 \,\,\mbox{\rm for}\,\, x \!\in\! (0,\frac{\pi}{4})
\Longrightarrow
\varphi(x)\!>\! 0 \,\,\mbox{\rm for}\,\, x \!\in\! (\frac{\pi}{4},\frac{\pi}{2}).
$

\bigskip
\noindent
{\boldmath $3^{0}$}
In this section of this appendix we give a proof based on the method from \cite{Malesevic_Makragic_2015} for the following equivalent inequality:
\begin{equation}
\label{Ineq_30}
\varphi(x)\!=\!\Big{(}\!\underbrace{-\frac{x^{6}}{\pi^{2}}+\frac{x^{4}}{2}-\frac{\pi^{2}}{16}x^{2}+\pi^{2}\!\Big{)}x}_{>0}\cos \!x+
\Big{(}\!\underbrace{\pi^{2}-4x^{2}}_{>0}\!\Big{)}(-\sin \!x)>0
\end{equation}
for $x \in (0,\frac{\pi}{2})$.

\noindent
Let us consider two cases:

\smallskip
\noindent
{\bf \boldmath $a$)}
If $x \!\in \!(0,c]$, where $c=\frac{3}{2}$, the following holds: $\cos x\geq\underline{T}_{\,6}^{\cos,0}(x)$ and $-\sin x\geq-\overline{T}_{5}^{\,\sin,0}(x)$, \cite{Malesevic_Makragic_2015}.
Then for $x\in \left(0,c\right]$ it holds:
\begin{equation}
\label{Ineq_31}
\begin{array}{rcl}
\varphi(x)
       \!&\!\!\!>\!\!\!&\!
\Big{(}\!\underbrace{-\frac{x^{6}}{\pi^{2}}+\frac{x^{4}}{2}-\frac{\pi^{2}}{16}x^{2}+\pi^{2}\Big{)}x}_{>0}\!\underline{T}_{\,6}^{\cos,0}(x)
                                                                                          \\[1.0 ex]
        \!&\!\!\!+\!\!\!&\!
{\Big(}\!\underbrace{\pi^{2}-4x^{2}}_{>0}\!{\Big)}\big{(}-\overline{T}_{5}^{\,\sin,0}(x)\big{)}
\!=\!P_{13}(x),
\end{array}
\end{equation}
where $P_{13}(x)$ is the polynomial
\begin{equation}
\label{Ineq_32}
\begin{array}{rcl}
P_{13}(x)
\!&\!\!=\!\!&\!
\mbox{\footnotesize $\displaystyle\frac{x^{3}}{11520\pi^{2}}$}
{\big (}
\mbox{\small $(16x^{10}\!+\!(\!-8\pi^{2}\!-\!480)x^{8}\!+\!(\pi^{4}\!+\!240\pi^{2}\!+\!5760)x^{6}$}
                                                                                              \\[2 ex]
\!&\!\!+\!\!&\!
\mbox{\small $(\!-46\pi^{4}\!-\!2496\pi^{2}\!-\!11520)x^{4}$}\!+\!\mbox{\small $(744\pi^{4}\!-\!1920\pi^{2})x^{2}$}
-\mbox{\small $\!4560\pi^{4}$}
                                                                                              \\[1 ex]
\!&\!\!+\!\!&\!
\mbox{\small $\!46080\,\pi^{2}\!$}
{\big )}\!=\!\mbox{\footnotesize $\displaystyle\frac{x^{3}}{11520\pi^{2}}$}
P_{10}(x).
\end{array}
\end{equation}
Then we determine the sign of the polynomial $P_{10}(x)$ for $x\in \left(0,c\right]$.
By introducing the substitute $z=x^{2}$, we get the fifth degree polynomial:
\begin{equation}
\begin{array}{rcl}
\label{Ineq_33}
P_{5}(z)
\!&\!\!\!\!=\!\!\!\!&\!
16z^{5}\!+\!(-8\pi^{2}-480)z^{4}
\!+\!(\pi^{4}\!+\!240\pi^{2}\!+\!5760)z^{3}
                                                                                           \\[2.0 ex]
\!&\!\!\!\!+\!\!\!\!&\!
\!(-46\pi^{4}\!\!-\!\!2496\pi^{2}\!\!-\!\!11520)z^{2}\!\!+\!\!(744\pi^{4}\!-\!1920\pi^{2})z
\!-\!4560\pi^{4}
                                                                                           \\[2.0 ex]
\!&\!\!\!\!+\!\!\!\!&\!
46080\pi^{2}
\,,
\end{array}
\end{equation}
and we determine the sign of the polynomial $P_{5}(z)$ for $z \in (0,c^{2}]$. Let us look at the first derivative of the polynomial $P_{5}(z)$, as the fourth degree polynomial, in the following form:
\begin{equation}
\begin{array}{rcl}
\label{Ineq_34}
P'_{5}(z)
\!&\!\!\!\!=\!\!\!\!&\!
80z^{4}\!+\!4(-8\pi^{2}-480)z^{3}
\!+\!3(\pi^{4}\!+\!240\pi^{2}\!+\!5760)z^{2}
                                                                                           \\[2.0 ex]
\!&\!\!\!\!+\!\!\!\!&\!
\!2(\!-46\pi^{4}\!\!-\!\!2496\pi^{2}\!\!-\!\!11520\!)z\!+\!744\pi^{4}\!-\!1920\pi^{2}\,.
\end{array}
\end{equation}
A real numerical factorization of the polynomial $P'_{5}(z)$, has been determined via {\sf Matlab} software, and given with $P'_{5}(z)\!=\!\alpha(z\!-\!z_{1})(z\!-\!z_{2})(z^{2}\!+\!pz\!+\!q)$,
where
$
\alpha\!=\!80,
z_{1}\!=\!0.871\ldots,
z_{2}\!=\!3.976\ldots,
$
$
p\!=\!-23.099\ldots,
q\!=\!193.019\ldots,
$
whereby the inequality \mbox{$p^{\,2}\!-4q\!<\!0$} is true.
The polynomial equation $P'_{5}(z)=0$ has got exactly two simple real roots with a symbolic radical representation and the corresponding numerical values $z_{1}$ and $z_{2}$. Let us notice that $z_{1}=0.871\ldots \in (0,c^{2}]$. Since $P_{5}'(0)=744\pi^{4}-1920\pi^{2}>0$ and $P'_{5}(1)=-7600-6224\pi^{2}+655\pi^{4}<0$
it follows that $P_{5}'(z)>0$ for $z<z_{1}$ and $P_{5}'(z)<0$ for $z \in (z_{1},z_{2})$.
Therefore, $P_{5}(z)$ is a monotonically increasing function for  $z \in (0,z_{1})$
and a monotonically decreasing function for $z \in (z_{1},c^{2}] \subset (z_{1},z_{2})$.
Hence, $P_{5}(z)$ reaches the maximum at the point $z_{1}\!=\!0.871\ldots$.
Then, since $P_{5}(0)=-4560\pi^{4}+46080\pi^{2}>0$ and
$P_{5}(c^{2})=-\frac{261711}{64}+\frac{1012887}{32}\pi^{2}-\frac{198879}{64}\pi^{4}>0$,
it follows that $P_{5}(z)$ is a positive function for $z \in (0,c^{2}]$.
Finally,
we conclude that
$P_{10}(x) \!>\! 0  \,\,\mbox{\rm for}\,\, x \!\in\! (0,c]
\Longrightarrow P_{13}(x) \!>\! 0 \,\,\mbox{\rm for}\,\, x \!\in\! (0,c]
\Longrightarrow \varphi(x)\!>\! 0 \,\,\mbox{\rm for}\,\, x \!\in\! (0,c]$.

\smallskip
\noindent
{\bf \boldmath $b$)}
If $x \!\in \!(c,\frac{\pi}{2})$, where $c=\frac{3}{2}$, let us define the function
\begin{equation}
\label{Ineq_35}
\begin{array}{rcl}
\phi(x)
          \!&\!\!\!=\!\!\!&\!
\varphi(\frac{\pi}{2}-x)=
\frac{(\pi-2x)}{8\pi^{2}}(\pi^{4}x^{2}+13\pi^{2}x^{4}+4x^{6})(-\sin x)
                                                                                          \\[1.5 ex]
           \!&\!\!\!+\!\!\!&\!
\frac{(\pi-2x)}{8\pi^{2}}(6\pi^{3}x^{3}+12\pi x^{5}+4\pi^{4})\sin x+(4\pi x-4x^{2})(-\cos x).
\end{array}
\end{equation}
Now we prove that $\phi(x)>0$ for $x \in (0,\hat{c})$, where
$\hat{c}=\frac{\pi}{2}-c=0.070\ldots$. The following holds:
$\!-\!\sin x\! \geq \!-\!\overline{T}_{1}^{\,\sin,0}(x)$,
$\sin x \! \geq \! \underline{T}_{\,3}^{\sin,0}(x)$ and
$\!-\!\cos x\! \geq\! -\!\overline{T}_{0}^{\,\cos,0}(x)$, \cite{Malesevic_Makragic_2015}.
Then for $x\in \left(0,\hat{c}\right)$ it holds:
\begin{equation}
\label{Ineq_36}
\begin{array}{rcl}
\phi(x)
       \!&\!\!\!\!>\!\!\!\!&\!
\frac{(\pi-2x)}{8\pi^{2}}(\pi^{4}x^{2}+13\pi^{2}x^{4}+4x^{6})
\big{(}\!-\!\overline{T}_{1}^{\,\sin,0}(x)\big{)}
                                                                                          \\[1.5 ex]
        \!&\!\!\!\!+\!\!\!\!&\!
\frac{(\pi\!-\!2x)}{8\pi^{2}}(6\pi^{3}x^{3}
\!+\!12\pi x^{5}\!+\!4\pi^{4})\underline{T}_{\,3}^{\sin,0}(\!x\!)
\!+\!(4\pi x\!-\!4x^{2})\big{(}\!\!-\!\!\overline{T}_{0}^{\,\cos,0}(\!x\!)\!\big{)}
                                                                                            \\[1.5 ex]
        \!&\!\!\!\!=\!\!\!\!&\!
Q_{9}(x),
\end{array}
\end{equation}
where $Q_{9}(x)$ is the polynomial
\begin{equation}
\label{Ineq_37}
\begin{array}{rcl}
Q_{9}(x)
          \!&\!\!\!\!=\!\!\!\!&\!
\mbox{\footnotesize $\displaystyle-\frac{x}{24\pi^{2}}$}
{\big (}\!
\mbox{\small $-12\pi x^{8}\!+\!(6\pi^{2}\!-\!24)x^{7}\!+\!(84\pi\!-\!6\pi^{3})x^{6}$}
\!+\!\mbox{\small $(3\pi^{4}\!-\!114\pi^{2})x^{5}$}
                                                                                                \\[1.5 ex]
           \!&\!\!\!\!+\!\!\!\!&\!
\mbox{\small $75\pi^{3}x^{4}$}\!-\!\mbox{\small $\!28\pi^{4}x^{3}\!+\!5\pi^{5}x^{2}$}
\!+\!\mbox{\small $(24\pi^{4}\!-\!96\pi^{2})x$}\!-\!
\mbox{\small $\!12\,\pi^{5}\!+\!96\,\pi^{3}\!$}
{\big )}
                                                                                                \\[1.5 ex]
           \!&\!\!\!\!=\!\!\!\!&\!
\mbox{\footnotesize $\displaystyle-\frac{x}{24\pi^{2}}$}Q_{8}(x).
\end{array}
\end{equation}
Then we determine the sign of the polynomial $Q_{8}(x)$ for $x\in \left(0,\hat{c}\right)$.
Let us look at the fourth derivative of the polynomial $Q_{8}(x)$, as the fourth degree polynomial, in the following form:
\begin{equation}
\begin{array}{rcl}
\label{Ineq_38}
Q^{(4)}_{8}(x)
\!&\!\!\!\!=\!\!\!\!&\!
-20160\pi x^{4}\!+\!840(6\pi^{2}-24)x^{3}
\!+\!360(84\pi\!-\!6\pi^{3})x^{2}
                                                                                           \\[2.0 ex]
\!&\!\!\!\!+\!\!\!\!&\!
\!120(3\pi^{4}\!\!-\!\!114\pi^{2})x\!+\!1800\pi^{3}\,.
\end{array}
\end{equation}
A real numerical factorization of the polynomial $Q^{(4)}_{8}(x)$, has been determined via {\sf Matlab} software, and given with $Q^{(4)}_{8}(x)\!=\!\beta(x\!-\!x_{1})(x\!-\!x_{2})(x^{2}\!+\!px\!+\!q)$,
where
$
\beta\!=\!-63334.507\ldots,
x_{1}\!=\!0.644\ldots,
x_{2}\!=\!-1.275\ldots,
$
$
p\!=\!-1.097\ldots,
q\!=\!1.071\ldots,
$
whereby the inequality \mbox{$p^{\,2}\!-4q\!<\!0$} is true.
The polynomial equation $Q^{(4)}_{8}(x)=0$ has got exactly two simple real roots with a symbolic radical representation and the corresponding numerical values $x_{1}$ and $x_{2}$.
Since $Q^{(4)}_{8}(0)=1800\pi^{3}>0$ it follows that $Q^{(4)}_{8}(x)$ is a positive function for
$x \in (0,\hat{c}) \subset (x_{2},x_{1})$.
Therefore, $Q^{(3)}_{8}(x)$ is a monotonically increasing function for $x \in (0,\hat{c})$.
Then, since $Q^{(3)}_{8}(\hat{c})<0$, it follows that $Q^{(3)}_{8}(x)$ is a negative function for
$x \in (0,\hat{c})$, so it follows that $Q^{(2)}_{8}(x)$ is a monotonically decreasing function for
$x \in (0,\hat{c})$. Since $Q^{(2)}_{8}(\hat{c})>0$, then it holds that $Q^{(2)}_{8}(x)$ is a positive function for $x \in (0,\hat{c})$. Hence, it follows that $Q'_{8}(x)$ is a monotonically increasing function for $x \in (0,\hat{c})$. Since $Q'_{8}(0)>0$, then it is valid that $Q'_{8}(x)$ is a positive function for $x \in (0,\hat{c})$. Hence, it holds that $Q_{8}(x)$ is a monotonically increasing function for $x \in (0,\hat{c})$. Finally, since
$Q_{8}(\hat{c})<0$ we conclude that
$Q_{8}(x) \!<\! 0  \,\,\mbox{\rm for}\,\, x \!\in\! (0,\hat{c})
\Longrightarrow Q_{9}(x) \!>\! 0 \,\,\mbox{\rm for}\,\, x \!\in\! (0,\hat{c})
\Longrightarrow \phi(x)\!>\! 0 \,\,\mbox{\rm for}\,\, x \!\in\! (0,\hat{c})
\Longrightarrow \varphi(x)\!>\! 0 \,\,\mbox{\rm for}\,\, x \!\in\! (c,\frac{\pi}{2}).
$

\bigskip
\noindent
{\boldmath $4^{0}$}
The inequality $f_{5}(x) \!>\! f_{4}(x)$, for $0 < x <\mbox{\small $\displaystyle\frac{\pi}{2}$}$,
is equivalent to the following trivial inequality
\begin{equation}
f_{5}(x) - f_{4}(x)
=
\displaystyle\frac{(\pi^2-4x^2)x^2}{16\pi^2} > 0\,.
\end{equation}

\bigskip

$\,$

\end{document}